\documentclass[a4paper,10pt,oneside,reqno]{amsart}

\usepackage[utf8]{inputenc}
\usepackage[explicit]{titlesec}
\usepackage{etoolbox}
\usepackage[hidelinks,bookmarksdepth=section]{hyperref}

\usepackage[a4paper]{geometry}
\usepackage{microtype}
\usepackage{amssymb,amsmath,amsopn,amsxtra,amsthm,amsfonts}
\usepackage{xr}
\usepackage[mathcal]{euscript}
\usepackage[bb=px]{mathalfa}

\titleformat{\section}
[hang] 
{\normalfont\scshape\filcenter} 
{\thesection.\ } 
{0em} 
{#1} 

\titleformat{\subsection}
[runin] 
{} 
{\bfseries\thesubsection.\ } 
{0em} 
{\itshape#1\ifstrempty{#1}{}{.}} 
[]

\titleformat{\subsubsection}
[runin] 
{} 
{\thesubsubsection.\ } 
{0em} 
{\itshape#1\ifstrempty{#1}{}{.}} 

\newtheorem{thm}[subsubsection]{Theorem}
\newtheorem{cor}[subsubsection]{Corollary}
\newtheorem{lem}[subsubsection]{Lemma}
\newtheorem{prop}[subsubsection]{Proposition}
\newtheorem{claim}[subsubsection]{Claim}

\theoremstyle{definition}

\theoremstyle{remark}
\newtheorem{rem}[subsubsection]{Remark}
\newtheorem{rems}[subsubsection]{Remarks}
\newtheorem{exam}[subsubsection]{Example}

\newcommand{\thmref}[1]{Theorem~\ref{#1}}
\newcommand{\secref}[1]{Sect.~\ref{#1}}
\newcommand{\ssecref}[1]{Subsect.~\ref{#1}}
\newcommand{\sssecref}[1]{\ref{#1}}
\newcommand{\lemref}[1]{Lemma~\ref{#1}}
\newcommand{\propref}[1]{Proposition~\ref{#1}}

\newcommand{\remref}[1]{Remark~\ref{#1}}

\renewcommand{\eqref}[1]{(\ref{#1})}

\newcommand{\examref}[1]{Example~\ref{#1}}

\newcommand{\itemref}[1]{\ref{#1}}

\numberwithin{equation}{subsection}

\usepackage{tikz}
\usetikzlibrary{matrix}
\usepackage{tikz-cd}

\usepackage{enumitem}
\setlist[enumerate,1]{label={(\alph*)},itemsep=\parskip,leftmargin=0pt}
\newlist{thmlist}{enumerate}{1}
\setlist[thmlist,1]{label={\em(\roman*)},ref={(\roman*)},itemsep=\parskip,leftmargin=0pt}     
\newlist{remlist}{enumerate}{1}
\setlist[remlist,1]{label={(\roman*)},itemsep=\parskip,leftmargin=0pt}

\usepackage{fixltx2e}

\usepackage{xspace}

\usepackage{xifthen}

\usepackage{xparse}

\usepackage{mathtools}

\newcommand{\nc}{\newcommand}
\nc{\renc}{\renewcommand}
\nc{\ssec}{\subsection}
\nc{\sssec}{\subsubsection}
\nc{\on}{\operatorname}
\nc{\term}[1]{#1\xspace}

\nc{\sA}{\ensuremath{\mathcal{A}}\xspace}
\nc{\sB}{\ensuremath{\mathcal{B}}\xspace}
\nc{\sC}{\ensuremath{\mathcal{C}}\xspace}
\nc{\sD}{\ensuremath{\mathcal{D}}\xspace}
\nc{\sE}{\ensuremath{\mathcal{E}}\xspace}
\nc{\sF}{\ensuremath{\mathcal{F}}\xspace}
\nc{\sG}{\ensuremath{\mathcal{G}}\xspace}
\nc{\sH}{\ensuremath{\mathcal{H}}\xspace}
\nc{\sI}{\ensuremath{\mathcal{I}}\xspace}
\nc{\sJ}{\ensuremath{\mathcal{J}}\xspace}
\nc{\sK}{\ensuremath{\mathcal{K}}\xspace}
\nc{\sL}{\ensuremath{\mathcal{L}}\xspace}
\nc{\sM}{\ensuremath{\mathcal{M}}\xspace}
\nc{\sN}{\ensuremath{\mathcal{N}}\xspace}
\nc{\sO}{\ensuremath{\mathcal{O}}\xspace}
\nc{\sP}{\ensuremath{\mathcal{P}}\xspace}
\nc{\sQ}{\ensuremath{\mathcal{Q}}\xspace}
\nc{\sR}{\ensuremath{\mathcal{R}}\xspace}
\nc{\sS}{\ensuremath{\mathcal{S}}\xspace}
\nc{\sT}{\ensuremath{\mathcal{T}}\xspace}
\nc{\sU}{\ensuremath{\mathcal{U}}\xspace}
\nc{\sV}{\ensuremath{\mathcal{V}}\xspace}
\nc{\sW}{\ensuremath{\mathcal{W}}\xspace}
\nc{\sX}{\ensuremath{\mathcal{X}}\xspace}
\nc{\sY}{\ensuremath{\mathcal{Y}}\xspace}
\nc{\sZ}{\ensuremath{\mathcal{Z}}\xspace}

\nc{\bA}{\ensuremath{\mathbf{A}}\xspace}
\nc{\bB}{\ensuremath{\mathbf{B}}\xspace}
\nc{\bC}{\ensuremath{\mathbf{C}}\xspace}
\nc{\bD}{\ensuremath{\mathbf{D}}\xspace}
\nc{\bE}{\ensuremath{\mathbf{E}}\xspace}
\nc{\bF}{\ensuremath{\mathbf{F}}\xspace}
\nc{\bG}{\ensuremath{\mathbf{G}}\xspace}
\nc{\bH}{\ensuremath{\mathbf{H}}\xspace}
\nc{\bI}{\ensuremath{\mathbf{I}}\xspace}
\nc{\bJ}{\ensuremath{\mathbf{J}}\xspace}
\nc{\bK}{\ensuremath{\mathbf{K}}\xspace}
\nc{\bL}{\ensuremath{\mathbf{L}}\xspace}
\nc{\bM}{\ensuremath{\mathbf{M}}\xspace}
\nc{\bN}{\ensuremath{\mathbf{N}}\xspace}
\nc{\bO}{\ensuremath{\mathbf{O}}\xspace}
\nc{\bP}{\ensuremath{\mathbf{P}}\xspace}
\nc{\bQ}{\ensuremath{\mathbf{Q}}\xspace}
\nc{\bR}{\ensuremath{\mathbf{R}}\xspace}
\nc{\bS}{\ensuremath{\mathbf{S}}\xspace}
\nc{\bT}{\ensuremath{\mathbf{T}}\xspace}
\nc{\bU}{\ensuremath{\mathbf{U}}\xspace}
\nc{\bV}{\ensuremath{\mathbf{V}}\xspace}
\nc{\bW}{\ensuremath{\mathbf{W}}\xspace}
\nc{\bX}{\ensuremath{\mathbf{X}}\xspace}
\nc{\bY}{\ensuremath{\mathbf{Y}}\xspace}
\nc{\bZ}{\ensuremath{\mathbf{Z}}\xspace}

\nc{\dA}{\ensuremath{\mathds{A}}\xspace}
\nc{\dB}{\ensuremath{\mathds{B}}\xspace}
\nc{\dC}{\ensuremath{\mathds{C}}\xspace}
\nc{\dD}{\ensuremath{\mathds{D}}\xspace}
\nc{\dE}{\ensuremath{\mathds{E}}\xspace}
\nc{\dF}{\ensuremath{\mathds{F}}\xspace}
\nc{\dG}{\ensuremath{\mathds{G}}\xspace}
\nc{\dH}{\ensuremath{\mathds{H}}\xspace}
\nc{\dI}{\ensuremath{\mathds{I}}\xspace}
\nc{\dJ}{\ensuremath{\mathds{J}}\xspace}
\nc{\dK}{\ensuremath{\mathds{K}}\xspace}
\nc{\dL}{\ensuremath{\mathds{L}}\xspace}
\nc{\dM}{\ensuremath{\mathds{M}}\xspace}
\nc{\dN}{\ensuremath{\mathds{N}}\xspace}
\nc{\dO}{\ensuremath{\mathds{O}}\xspace}
\nc{\dP}{\ensuremath{\mathds{P}}\xspace}
\nc{\dQ}{\ensuremath{\mathds{Q}}\xspace}
\nc{\dR}{\ensuremath{\mathds{R}}\xspace}
\nc{\dS}{\ensuremath{\mathds{S}}\xspace}
\nc{\dT}{\ensuremath{\mathds{T}}\xspace}
\nc{\dU}{\ensuremath{\mathds{U}}\xspace}
\nc{\dV}{\ensuremath{\mathds{V}}\xspace}
\nc{\dW}{\ensuremath{\mathds{W}}\xspace}
\nc{\dX}{\ensuremath{\mathds{X}}\xspace}
\nc{\dY}{\ensuremath{\mathds{Y}}\xspace}
\nc{\dZ}{\ensuremath{\mathds{Z}}\xspace}

\nc{\bbA}{\ensuremath{\mathbb{A}}\xspace}
\nc{\bbB}{\ensuremath{\mathbb{B}}\xspace}
\nc{\bbC}{\ensuremath{\mathbb{C}}\xspace}
\nc{\bbD}{\ensuremath{\mathbb{D}}\xspace}
\nc{\bbE}{\ensuremath{\mathbb{E}}\xspace}
\nc{\bbF}{\ensuremath{\mathbb{F}}\xspace}
\nc{\bbG}{\ensuremath{\mathbb{G}}\xspace}
\nc{\bbH}{\ensuremath{\mathbb{H}}\xspace}
\nc{\bbI}{\ensuremath{\mathbb{I}}\xspace}
\nc{\bbJ}{\ensuremath{\mathbb{J}}\xspace}
\nc{\bbK}{\ensuremath{\mathbb{K}}\xspace}
\nc{\bbL}{\ensuremath{\mathbb{L}}\xspace}
\nc{\bbM}{\ensuremath{\mathbb{M}}\xspace}
\nc{\bbN}{\ensuremath{\mathbb{N}}\xspace}
\nc{\bbO}{\ensuremath{\mathbb{O}}\xspace}
\nc{\bbP}{\ensuremath{\mathbb{P}}\xspace}
\nc{\bbQ}{\ensuremath{\mathbb{Q}}\xspace}
\nc{\bbR}{\ensuremath{\mathbb{R}}\xspace}
\nc{\bbS}{\ensuremath{\mathbb{S}}\xspace}
\nc{\bbT}{\ensuremath{\mathbb{T}}\xspace}
\nc{\bbU}{\ensuremath{\mathbb{U}}\xspace}
\nc{\bbV}{\ensuremath{\mathbb{V}}\xspace}
\nc{\bbW}{\ensuremath{\mathbb{W}}\xspace}
\nc{\bbX}{\ensuremath{\mathbb{X}}\xspace}
\nc{\bbY}{\ensuremath{\mathbb{Y}}\xspace}
\nc{\bbZ}{\ensuremath{\mathbb{Z}}\xspace}

\nc{\mrm}[1]{\ensuremath{\mathrm{#1}}\xspace}
\nc{\mit}[1]{\ensuremath{\mathit{#1}}\xspace}
\nc{\mbf}[1]{\ensuremath{\mathbf{#1}}\xspace}
\nc{\mcal}[1]{\ensuremath{\mathcal{#1}}\xspace}
\nc{\msc}[1]{\ensuremath{\mathscr{#1}}\xspace}

\renc{\bar}[1]{\overline{#1}}

\let\sectsign\S
\let\S\relax

\nc{\sub}{\subset}
\nc{\too}{\longrightarrow}
\nc{\hook}{\hookrightarrow}
\nc*{\hooklongrightarrow}{\ensuremath{\lhook\joinrel\relbar\joinrel\rightarrow}}
\nc{\hooklong}{\hooklongrightarrow}
\nc{\twoheadlongrightarrow}{\relbar\joinrel\twoheadrightarrow}
\nc{\isoto}{\xrightarrow{\sim}}
\nc{\isofrom}{\xleftarrow{\sim}}
\renc{\ge}{\geqslant}
\renc{\geq}{\geqslant}
\renc{\le}{\leqslant}
\renc{\leq}{\leqslant}

\nc{\id}{\mathrm{id}}

\DeclareMathOperator{\Hom}{\on{Hom}}
\nc{\uHom}{\underline{\smash{\Hom}}}
\DeclareMathOperator{\Maps}{\on{Maps}}

\DeclareMathOperator{\End}{\on{End}}
\DeclareMathOperator{\Sym}{\on{Sym}}
\nc{\uEnd}{\underline{\smash{\End}}}
\DeclareMathOperator{\codim}{\on{codim}}
\nc{\colim}{\varinjlim}
\renc{\lim}{\varprojlim}
\nc{\Cofib}{\on{Cofib}}
\nc{\Fib}{\on{Fib}}
\nc{\initial}{\varnothing}
\nc{\op}{\mathrm{op}}

\DeclareMathOperator*{\fibprod}{\times}

\renc{\setminus}{\smallsetminus}

\newcommand{\spref}[1]{\href{http://stacks.math.columbia.edu/tag/#1}{#1}}
\newcommand{\spcite}[1]{\cite[\spref{#1}]{stacks-project}}

\nc{\CRing}{\mrm{CRing}}
\nc{\Poly}{\mrm{Poly}}
\nc{\Spc}{\mrm{Spc}}
\nc{\pt}{\mrm{pt}}
\nc{\SCRing}{\mrm{SCRing}}
\nc{\Mod}{\mrm{Mod}}
\nc{\Spt}{\mrm{Spt}}
\nc{\A}{\bA}
\nc{\red}{{\mrm{red}}}
\nc{\Spec}{\on{Spec}}
\nc{\Coh}{\on{Coh}}
\nc{\cn}{{\on{cn}}}
\nc{\Sm}{\mrm{Sm}}
\nc{\Sch}{\mrm{Sch}}
\nc{\aff}{{\mrm{aff}}}
\nc{\Einfty}{\sE_\infty}
\nc{\Funct}{\on{Funct}}
\nc{\sift}{{\mrm{sift}}}
\nc{\dash}{{\textnormal{-}}}
\nc{\dashalg}{\dash\mrm{alg}}
\nc{\EInfAlgCn}{\Einfty\dashalg^{\cn}}
\nc{\SCRingMod}{\mrm{SCRing\textnormal{+}Mod}}
\nc{\DSch}{\mrm{DSch}}
\nc{\Qcoh}{\on{Qcoh}}
\nc{\InftyCat}{\infty\dash\mrm{Cat}}
\nc{\DStk}{\mrm{DStk}}
\nc{\open}{\mrm{open}}
\nc{\bDelta}{\mathbf{\Delta}}
\nc{\Cech}{\textnormal{\v{C}}}
\nc{\Perf}{\on{Perf}}
\nc{\cl}{{\mrm{cl}}}
\nc{\Tot}{\on{Tot}}
\nc{\Sh}{\on{Sh}}
\nc{\comp}{{\mrm{comp}}}
\nc{\idem}{{\mrm{idem}}}
\nc{\Ind}{{\on{Ind}}}
\nc{\Open}{\mrm{Open}}
\nc{\Wald}{\mrm{Wald}}
\nc{\Stab}{\mrm{Stab}}
\nc{\iso}{\on{iso}}
\nc{\cofib}{\mrm{cofib}}
\nc{\K}{\on{K}}
\nc{\qcqs}{{\mrm{qcqs}}}
\nc{\coh}{{\mrm{coh}}}
\nc{\perf}{{\mrm{perf}}}
\nc{\Der}{\on{Der}}
\nc{\triv}{\mrm{triv}}
\nc{\obstr}{\mrm{obstr}}
\nc{\pr}{\mrm{pr}}
\renc{\P}{\bP}
\nc{\Proj}{{\on{Proj}}} 
\nc{\QcohAlg}{\on{QcohAlg}}
\nc{\modmod}{/\!\!/}
\nc{\Bl}{\mrm{Bl}}
\nc{\Res}{\mrm{Res}}
\nc{\univ}{\mrm{univ}}
\nc{\closed}{{\mrm{closed}}}
\nc{\lfp}{{\mrm{lfp}}}
\nc{\M}{\sM}
\nc{\Mbar}{\bar{\M}}
\nc{\DM}{{\mrm{DM}}}
\nc{\lfta}{{\mrm{lfta}}}
\nc{\proj}{{\mrm{proj}}}
\nc{\Pro}{\on{Pro}}
\nc{\gp}{{\mrm{gp}}}
\nc{\laxprod}{{\rightarrow}}
\nc{\Tor}{\on{Tor}}
\nc{\der}{{\mrm{der}}}
\nc{\Kosz}{\on{Kosz}}
\nc{\Ann}{\on{Ann}}
\nc{\codimvir}{\on{codim.\!vir}}
\nc{\codimtop}{\on{codim.\!top}}
\nc{\Div}{\mrm{Div}}
\nc{\GDiv}{\mrm{GDiv}}
\nc{\VDiv}{\mrm{VDiv}}
\nc{\uPic}{\underline{\mrm{Pic}}}

\nc{\scr}{\term{simplicial commutative ring}}
\nc{\scrs}{\term{simplicial commutative rings}}

\nc{\ccrs}{\term{connective $\Einfty$-ring spectrum}}
\nc{\ccrss}{\term{connective $\Einfty$-ring spectra}}

\nc{\inftyCat}{\term{$\infty$-category}}
\nc{\inftyCats}{\term{$\infty$-categories}}

\nc{\inftyGrpd}{\term{$\infty$-groupoid}}
\nc{\inftyGrpds}{\term{$\infty$-groupoids}}

\title{Virtual~Cartier~divisors~and~blow-ups}
\author{Adeel A. Khan}
\address{Fakultät für Mathematik\\
Universität Regensburg\\
93040 Regensburg\\
Germany}
\email{\href{mailto:adeel.khan@mathematik.uni-regensburg.de}{adeel.khan@mathematik.uni-regensburg.de}}

\author{David Rydh}
\address{KTH Royal Institute of Technology\\
Department of Mathematics\\
SE-100 44 Stockholm\\
Sweden}
\email{\href{mailto:dary@math.kth.se}{dary@math.kth.se}}

\date{May 16, 2019}

\begin{document}


\maketitle

\parskip 0pt
\parskip 0.2cm


\section{Introduction}
\label{sec:intro}

\ssec{}

Let $X$ be a scheme, $Z$ a closed subscheme, and $\tilde{X}$ the blow-up of $X$ in $Z$.
Recall that $\tilde{X}$ admits the following universal property: for any morphism $f : S \to X$ such that the schematic fibre $f^{-1}(Z)$ is an effective Cartier divisor on $S$, there exists a unique morphism $S \to \tilde{X}$ over $X$.

The purpose of this note is to prove a stronger universal property, where the morphism $f$ is allowed to be arbitrary, when $Z\to X$ is a regular immersion.
In particular, we provide a complete description of the functor represented by the blow-up $\tilde{X}$.
Namely, we show that there is a canonical bijection between the set of $X$-morphisms $S \to \tilde{X}$ and the set of \emph{virtual effective Cartier divisors on $S$ lying over $(X,Z)$}.
In a word, a virtual effective Cartier divisor is a closed subscheme that is equipped with some additional structure that remembers that it is cut out locally by a single equation (or ``of virtual codimension $1$'').

\ssec{}

We also use the notion of virtual Cartier divisors to construct blow-ups of quasi-smooth closed immersions of \emph{derived} schemes and stacks.
This generalizes a local construction that was used by Kerz--Strunk--Tamme in their proof of Weibel's conjecture on negative K-theory \cite{KST} (see \sssecref{sssec:KST}).
In \cite{KhanKBlow} quasi-smooth blow-ups are used to give a simple new proof of Cisinski's cdh descent for homotopy K-theory.
Quasi-smooth blow-ups also give rise to a derived version of Verdier's technique of deformation to the normal cone (\thmref{thm:deformation}).
Using this one can give a simple new construction of virtual fundamental classes that does not use the intrinsic normal cone of Behrend--Fantechi.

\ssec{}

The organization of this paper is as follows.
In \secref{sec:regular} we study regular closed immersions from the perspective of derived algebraic geometry; this material is well-known to the experts.
In \secref{sec:divisors} we define virtual effective Cartier divisors and show that they coincide with generalized effective Cartier divisors.
\secref{sec:blow-up} contains our results on blow-ups.
The derived blow-up is constructed in
\sssecref{sssec:blow-up}/\sssecref{sssec:Bl_Z/X} and its properties are summarized in \thmref{thm:all}.
The universal property mentioned above is in \sssecref{cor:universal property in classical case}.
The remainder of the section is concerned with the proofs.
Finally, in \secref{sec:multiple} we generalize the construction to simultaneous blow-ups in multiple centres.
In \ssecref{ssec:loc-reg-emb} we discuss local regular immersions, that is, quasi-smooth finite unramified morphisms, which generalize regular closed immersions.
The main properties of simultaneous blow-ups are summarized in \thmref{thm:multiple}.

\ssec{}

We would like to thank Denis-Charles Cisinski and Akhil Mathew for comments on previous revisions.


\section{Quasi-smooth immersions}
\label{sec:regular}

\ssec{}
\label{ssec:homological}

We begin by reviewing the notion of regular closed immersion in classical algebraic geometry.

\sssec{}
Let $A$ be a commutative ring.
For an element $f \in A$, the Koszul complex $\Kosz_A(f)$ is the chain complex
  \begin{equation*}
    \Kosz_A(f) := \left(A \xrightarrow{f} A\right),
  \end{equation*}
concentrated in degrees $0$ and $1$.
Thus $H_0(\Kosz_A(f)) = A/f$ and $H_1(\Kosz_A(f)) = \Ann_A(f)$ is the annihilator.
In particular $\Kosz_A(f)$ is acyclic in positive degrees, and hence quasi-isomorphic to $A/f$, if and only if $f$ is regular (a non-zero divisor).
More generally, given a sequence of elements $(f_1,\ldots,f_n)$, the Koszul complex $\Kosz_A(f_1,\ldots,f_n)$ is defined as the tensor product (over $A$)
  \begin{equation*}
    \Kosz_A(f_1,\ldots,f_n) = \bigotimes_{i} \left(A \xrightarrow{f_i} A\right).
  \end{equation*}
We say that the sequence $(f_1,\ldots,f_n)$ is \emph{regular} if the Koszul complex is acyclic in positive degrees.
This is called a Koszul-regular sequence in the Stacks Project~\spcite{062D}.
When $A$ is noetherian and $f_i$ belong to the radical, this is equivalent to the usual inductive definition: $f_1$ is regular, $f_2$ is regular in $A/(f_1)$, etc. (see \cite[Cor.~19.5.2]{EGAIV4} and \cite[Prop.~1.3]{SGA6}).

\sssec{}
Any sequence $(f_1,\ldots,f_n)$ determines a homomorphism $\bZ[T_1,\ldots,T_n] \to A$, $T_i \mapsto f_i$, and the Koszul complex $\Kosz_A(f_1,\ldots,f_n)$ is quasi-isomorphic to the derived tensor product
  \begin{equation*}
    A \otimes^{\bL}_{\bZ[T_1,\ldots,T_n]} \bZ[T_1,\ldots,T_n]/(T_1,\ldots,T_n);
  \end{equation*}
indeed the sequence $(T_1,\ldots,T_n)$ is regular, so that $\Kosz_{\bZ[T_1,\ldots,T_n]}(T_1,\ldots,T_n)$ provides a free resolution of $\bZ[T_1,\ldots,T_n]/(T_1,\ldots,T_n)$, and the Koszul complex is stable under arbitrary extension of scalars.

\sssec{}
Let $i : Z\hook X$ be a closed immersion of schemes.
We say that $i$ is \emph{regular} if its ideal of definition $\sI \subset \sO_X$ is Zariski-locally generated by a regular sequence.
This is equivalent to the definition in \cite[Exp.~VII, D\'ef.~1.4]{SGA6}, and to the definition in \cite[D\'ef.~16.9.2]{EGAIV4} when $X$ is locally noetherian.
In the Stacks Project this is called a Koszul-regular ideal~\spcite{07CV}.
When $i$ is regular, the conormal sheaf $\sN_{Z/X} = \sI/\sI^2$ is locally free of finite rank.
Moreover, the relative cotangent complex $\sL_{Z/X}$ is canonically identified with $\sN_{Z/X}[1]$.
For us an \emph{effective Cartier divisor} on a scheme $X$ will be a scheme $D$ equipped with a regular closed immersion $i_D : D \hook X$ of codimension $1$.

\ssec{}

We now re-interpret the above discussion in the language of derived algebraic geometry.
The basic idea is that the Koszul complex $\Kosz_A(f_1,\ldots,f_n)$ can be viewed as the ``ring of functions'' on a certain \emph{derived} subscheme of $\Spec(A)$.
In order to make sense of this, one should work with simplicial commutative rings, which following Quillen is the natural setting for derived tensor products of commutative algebras.

\sssec{}
We use the language of \inftyCats; our reference is \cite{HTT}.
Let $\Spc$ be the \inftyCat of spaces\footnotemark.
\footnotetext{We will not rely on any particular model, say by topological spaces or simplicial sets; instead we will use the language of \inftyGrpds.
Thus a ``space'' has ``points'' (objects), ``paths'' (invertible morphisms) between any two points, etc.
In particular, the term ``isomorphism of spaces'' will never refer to an isomorphism of set-theoretic models, but rather to an isomorphism in the \inftyCat $\Spc$.}
For any \inftyCat $\bC$, there are mapping spaces $\Maps_\bC(x,y) \in \Spc$ for any pair of objects $x,y \in \bC$.

\sssec{}
Let $\SCRing$ be the \inftyCat of \scrs\footnotemark;
\footnotetext{This is the \inftyCat of functors $(\Poly)^\op \to \Spc$, from the opposite of the category $\Poly$ of polynomial algebras $\bZ[T_1,\ldots,T_n]$ ($n\ge 0$) and ring homomorphisms, to the \inftyCat of spaces, which send finite coproducts in $\Poly$ to finite products of spaces.
As for spaces, we will not use any set-theoretic models; the term ``simplicial commutative ring'' just means ``object of $\SCRing$''.
We can think of such an object $R : (\Poly)^\op \to \Spc$ as an underlying space $R(\bZ[T]) \in \Spc$, equipped with certain operations as encoded by the category $\Poly$.}
we refer to \cite[\sectsign~25.1]{SAG-20180204} for a detailed account.
If we forget the multiplication on a \scr $A$, we get by the Dold--Kan correspondence an ``underlying chain complex''.
If we forget both the addition and multiplication we get an ``underlying space''.
We say that a \scr $A$ is $0$-truncated or ``discrete'' if the underlying space (pointed at $0$) has no higher homotopy groups (i.e., can be identified with a set); the full subcategory of discrete simplicial commutative rings is canonically equivalent to the category $\CRing$ of ordinary commutative rings, and under this identification the $0$-truncation functor is given by $A \mapsto \pi_0(A)$.

Using Quillen's machinery of non-abelian derived functors one defines derived tensor products of \scrs: they are computed using simplicial resolutions by polynomial algebras $\bZ[T_1,\ldots,T_n]$.
For example, the derived tensor product
  \begin{equation*}
    A \otimes^{\bL}_{\bZ[T_1,\ldots,T_n]} \bZ[T_1,\ldots,T_n]/(T_1,\ldots,T_n)
  \end{equation*}
can be viewed as a simplicial commutative ring, for any commutative ring $A$ and any homomorphism $\bZ[T_1,\ldots,T_n] \to A$.
The underlying chain complex recovers the Koszul complex of the corresponding elements (as an object of $\Mod_\bZ$, the stable \inftyCat of chain complexes of $\bZ$-modules).

\sssec{}
Derived algebraic geometry is an extension of classical algebraic geometry that allows one to make sense of ``$\Spec(A)$'' where $A$ is a simplicial commutative ring.
We refer to \cite{SAG-20180204} or \cite{HAG2} for details.

The object $\Spec(A)$ is an affine derived scheme which has an ``underlying classical scheme'' $\Spec(A)_\cl = \Spec(\pi_0(A))$, but is further equipped with a quasi-coherent sheaf of \scrs $\sO_{\Spec(A)}$.
Thus it can be viewed as an infinitesimal thickening of $\Spec(A)_\cl$ where the nilpotents live in the higher homotopy groups.
For example, given affine schemes $X = \Spec(A)$ and $Y= \Spec(B)$ over $S = \Spec(R)$, there is a derived version of the fibred product:
  \begin{equation*}
    X \fibprod^{\bR}_S Y = \Spec(A \otimes^\bL_R B),
  \end{equation*}
which is a derived scheme with $(X\fibprod^{\bR}_S Y)_\cl = X \fibprod_S Y$.

A general derived scheme is Zariski-locally of the form $\Spec(A)$ for $A \in \SCRing$.
Derived fibred products are homotopy limits in the \inftyCat of derived schemes.
The category of classical schemes embeds fully faithfully into the \inftyCat of derived schemes, and its essential image is spanned by derived schemes $X$ whose structure sheaf $\sO_X$ is discrete (i.e., takes values in discrete \scrs).
This embedding is left adjoint to $(-)\to (-)_\cl$ and hence preserves colimits.
The embedding does not preserve fibred products though.
We will say that a derived scheme is ``classical'' if it is in the essential image.

The discussion of \ssecref{ssec:homological} can now be rephrased as follows.

\begin{prop}\label{prop:regular immersion of classical schemes}
Let $i : Z \hook X$ be a closed immersion of schemes.
Then $i$ is regular if and only if Zariski-locally on $X$, there exists a morphism $f : X \to \A^n$ and a commutative square
  \begin{equation*}
    \begin{tikzcd}
      Z \ar[hookrightarrow]{r}\ar{d}
        & X\ar{d}
      \\
      \{0\}\ar[hookrightarrow]{r}
        & \A^n,
    \end{tikzcd}
  \end{equation*}
which is homotopy cartesian in the \inftyCat of derived schemes.
\end{prop}

Here $\A^n = \Spec(\bZ[T_1,\ldots,T_n])$ denotes $n$-dimensional affine space over $\Spec(\bZ)$, and $\{0\} = \Spec(\bZ[T_1,\ldots,T_n]/(T_1,\ldots,T_n))$ denotes the inclusion of the origin.

\ssec{}

We will now extend the notion of regularity to the derived setting.

\sssec{}
Let $A$ be a \scr.
Let $f_1,\ldots,f_n \in A$ be a sequence of elements (i.e., points in the underlying space of $A$).
Let $A\modmod(f_1,\ldots,f_n)$ denote the \scr defined by the homotopy cocartesian square
  \begin{equation}\label{eq:A//f}
    \begin{tikzcd}
      \bZ[T_1,\ldots,T_n]\ar{r}\ar{d}{T_i\mapsto f_i}
        & \bZ[T_1,\ldots,T_n]/(T_1,\ldots,T_n) \ar{d}
      \\
      A \ar{r}
        & A\modmod(f_1,\ldots,f_n).
    \end{tikzcd}
  \end{equation}
That is, $A\modmod(f_1,\ldots,f_n)$ is given by the derived tensor product
\begin{equation*}
A \otimes^\bL_{\bZ[T_1,\ldots,T_n]} \bZ[T_1,\ldots,T_n]/(T_1,\ldots,T_n).
\end{equation*}
We have $\pi_0(A\modmod(f_1,\ldots,f_n)) = \pi_0(A)/(f_1,\ldots,f_n)$ where we write $f_i$ again for its connected component in $\pi_0(A)$.
The underlying $A$-module of $A\modmod(f_1,\ldots,f_n)$ is given by
  \begin{equation*}
    \bigotimes^\bL_i \Cofib(A\xrightarrow{f_i} A),
  \end{equation*}
where $\Cofib$ denotes the homotopy cofibre (in the stable \inftyCat of $A$-modules).

\begin{exam}
If $A$ is discrete and the sequence $(f_1,\ldots,f_n)$ is regular, then the canonical homomorphism $A\modmod(f_1,\ldots,f_n) \to A/(f_1,\ldots,f_n)$ is an isomorphism.
\end{exam}

\begin{exam}
Let $A$ be discrete and consider $A\modmod(0) \in \SCRing$.
Its underlying chain complex is given by $A\oplus A[1]$ (with zero differential).
In particular we have $\pi_0(A\modmod(0)) = \pi_1(A\modmod(0)) = A$.
\end{exam}

\sssec{}
By construction, $A\modmod(f_1,\ldots,f_n)$ admits the following universal property:

\begin{lem} \label{lem:map-from-A//f}
Let $A$ be a \scr and $f_1,\ldots,f_n \in A$ a sequence of points in the underlying space of $A$.
Then for any $A$-algebra $B \in \SCRing_A$, there is a canonical isomorphism of spaces
  \begin{equation*}
    \Maps_{\SCRing_A}(A\modmod(f_1,\ldots,f_n), B)
      \simeq \prod_{1\le i\le n} \Maps_{B}(f'_i, 0),
  \end{equation*}
where $f'_i$ are the images of $f_i$ in $B$.
That is, the space of $A$-algebra homomorphisms $A\modmod(f_1,\ldots,f_n) \to B$ is isomorphic to the space of paths $f'_i\simeq 0$, for each $1\le i\le n$, in the underlying space of $B$.
\end{lem}

\begin{proof}
For any $A$-algebra $B$, the space of commutative squares
  \begin{equation*}
    \begin{tikzcd}
      \bZ[T_1,\ldots,T_n]\ar{r}\ar{d}{T_i\mapsto f_i}
        & \bZ[T_1,\ldots,T_n]/(T_1,\ldots,T_n) \ar{d}
      \\
      A \ar{r}
        & B
    \end{tikzcd}
  \end{equation*}
is nothing else than the space of identifications between the two possible composites $\bZ[T_1,\ldots,T_n] \to B$.
Since $\bZ[T_1,\ldots,T_n]$ is free as a \scr, this is the same as the space of identifications between the two sequences $(f_1,\ldots,f_n)$ and $(0,\ldots,0)$ of points in the underlying space of $B$.
\end{proof}

\sssec{}
Let $i : Z \hook X$ be a closed immersion of derived schemes (i.e., the underlying morphism of classical schemes $i_\cl : Z_\cl \to X_\cl$ is a closed immersion).
We say that $i$ is \emph{quasi-smooth} if Zariski-locally on $X$, there exists a morphism $f : X \to \A^n$ and a homotopy cartesian square
  \begin{equation*}
    \begin{tikzcd}
      Z \ar[hookrightarrow]{r}\ar{d}
        & X\ar{d}
      \\
      \{0\}\ar[hookrightarrow]{r}
        & \A^n,
    \end{tikzcd}
  \end{equation*}
in the \inftyCat of derived schemes.
In other words, $Z \hook X$ is locally of the form $\Spec(A\modmod(f_1,\ldots,f_n)) \hook \Spec(A)$, for some $f_1,\ldots,f_n \in A$.
If $X$ and $Z$ are classical, then a closed immersion $i : Z \hook X$ is quasi-smooth if and only if it is a regular immersion.
By \propref{prop:regular immersion of classical schemes} this agrees with our previous definition when $X$ and $Z$ are classical.
However, even if $X$ is classical, there exist quasi-smooth immersions $Z \hook X$ with $Z$ non-classical (as long as $X$ is nonempty).

\sssec{}
We now give a differential characterization of quasi-smooth immersions.

\begin{prop}\label{prop:regular immersions and cotangent}
Let $i : Z \hook X$ be a closed immersion of derived schemes.
Then $i$ is quasi-smooth if and only if it is locally of finite presentation and the shifted cotangent complex $\sL_{Z/X}[-1]$ is a locally free $\sO_Z$-module of finite rank.
\end{prop}

\begin{proof}
The condition is clearly necessary: since it is Zariski-local and stable under arbitrary derived base change, we may assume that $i$ is the inclusion of the origin $\{0\} \hook \A^n$, $n\ge 0$, in which case $\sL_{\{0\}/\A^n}[-1] = \sN_{\{0\}/\A^n}$ is free of rank $n$.

Conversely, suppose that $X = \Spec(A)$ and $Z = \Spec(B)$ are affine, and the shifted cotangent complex $\mrm{L}_{B/A}[-1] \in \Mod_B$ is free of rank $n$.
Let $F$ denote the homotopy fibre of the morphism $\varphi : A \to B$ in the stable \inftyCat $\Mod_A$, so that there is a canonical isomorphism of $\pi_0(B)$-modules $\pi_1(\mrm{L}_{B/A})\simeq \pi_0(F \otimes^\bL_A B)$ \cite[Cor.~25.3.6.1]{SAG-20180204}.
Choose a basis $df_1,\ldots,df_n$ for $\pi_1(\mrm{L}_{B/A})$ and note that the corresponding elements of $\pi_0(F \otimes^\bL_A B)$ lift to elements $\tilde{f}_1,\ldots,\tilde{f}_n \in \pi_0(F)$, since $\varphi : A \to B$ is surjective on $\pi_0$; moreover, we can assume by Nakayama's lemma that the $\tilde{f_i}$ generate $\pi_0(F)$ as a $\pi_0(A)$-module.
Lifting them to points in the underlying space of $F$, we get points $f_i \in A$ equipped with paths $\varphi(f_i) \simeq 0$ in $B$, and hence a canonical homomorphism of \scrs $A\modmod(f_1,\ldots,f_n) \to B$ (\lemref{lem:map-from-A//f}).
By construction it induces an isomorphism $\pi_0(A)/(f_1,\ldots,f_n) \simeq \pi_0(B)$ on connected components, so by \cite[Cor.~25.3.6.6]{SAG-20180204} it suffices to show that its relative cotangent complex vanishes, which follows by examining the exact triangle
  \begin{equation*}
    \mrm{L}_{(A\modmod(f_i)_i)/A} \otimes^\bL_{A\modmod(f_i)_i} B \to \mrm{L}_{B/A} \to \mrm{L}_{B/(A\modmod(f_i)_i)}
  \end{equation*}
in $\Mod_B$.
\end{proof}

Given a quasi-smooth closed immersion $i : Z \hook X$, we write $\sN_{Z/X} = \sL_{Z/X}[-1]$ and take this as the definition of the \emph{conormal sheaf}.
By the above characterization, this is a locally free $\sO_Z$-module of finite rank.
It also follows that quasi-smoothness is an fpqc-local property, in view of \cite[Lem.~2.9.1.4]{SAG-20180204}.

\begin{exam}\label{ex:between-smooth-is-qs}
If $X$ and $Z$ are \emph{smooth} over some base $S$, then any closed immersion $i : Z \hook X$ is quasi-smooth.
This follows from the exact triangle
  \begin{equation*}
    i^*\sL_{X/S} \to \sL_{Z/S} \to \sL_{Z/X}.
  \end{equation*}
\end{exam}

\begin{exam}
If $i$ admits a smooth retraction, then it is quasi-smooth.
This is a special case of \examref{ex:between-smooth-is-qs}.
\end{exam}

\sssec{}
Let $i : Z \hook X$ be a quasi-smooth closed immersion of derived schemes.
The \emph{virtual codimension} of $i$, defined Zariski-locally on $Z$, is the rank of the locally free $\sO_Z$-module $\sN_{Z/X}$.
We will sometimes denote it by $\codimvir(Z,X)$.
Locally this number is determined by the formula $\codimvir(\Spec(A\modmod(f_1,\ldots,f_n)), \Spec(A)) = n$ for any $A \in \SCRing$ and points $f_1,\ldots,f_n \in A$.
Note that virtual codimension is stable under arbitrary derived base change: $\codimvir(Z,X) = \codimvir(Z \fibprod^\bR_X X', X')$ for any morphism $f : X' \to X$ of derived schemes.

\sssec{}
Define the \emph{topological (Krull) codimension} $\codimtop(Z,X)$ as the topological codimension of $Z_\cl$ in $X_\cl$.
If $X_\cl$ is locally noetherian then we have an inequality $\codimvir(Z,X) \ge \codimtop(Z,X)$.
This is an equality (at a point~$x \in X$) if the derived scheme $Z \fibprod^\bR_X X_\cl$ is classical (in a Zariski neighbourhood of~$x$).
If $X_\cl$ is Cohen--Macaulay (e.g.\ regular) at the point $x$, then the converse also holds.

\sssec{}

A morphism of derived schemes $f: Y \to X$ is called \emph{quasi-smooth} if it admits, Zariski-locally on $Y$, a factorization
  \begin{equation*}
    Y \xrightarrow{i} X' \xrightarrow{p} X,
  \end{equation*}
with $i$ a quasi-smooth closed immersion and $p$ smooth.
Quasi-smooth morphisms between classical schemes are usually called lci (local complete intersection).

\begin{prop}
Let $f : Y \to X$ be a morphism of derived schemes.
Then $f$ is quasi-smooth if and only if it is locally of finite presentation and the cotangent complex $\sL_{Y/X}$ is of Tor-amplitude $\le 1$.
\end{prop}

\begin{proof}
By \cite[Prop.~2.8.4.2]{SAG-20180204}, the question is local on $Y$.
If $f$ admits a factorization as above, then it is clear that $f$ is locally of finite presentation.
The exact triangle $i^*\sL_{X'/X} \to \sL_{Y/X} \to \sL_{Y/X'}$ shows that $\sL_{Y/X}$ is also of Tor-amplitude $\le 1$.
For the converse direction, let $\varphi : A \to B$ be a homomorphism of \scrs which is locally of finite presentation and such that $L_{B/A}$ is of Tor-amplitude $\le 1$.
Then $\pi_0(A) \to \pi_0(B)$ is of finite presentation in the sense of ordinary commutative algebra, so we may find a homomorphism $A' = A[T_1,\ldots,T_n] \to B$, for some $n \ge 0$, that extends $\varphi$ and is surjective on $\pi_0$.
The homomorphism $A \to A'$ is smooth, so it will suffice to show that $L_{B/A'}[-1]$ is locally free of finite rank (or equivalently, since $\pi_0(L_{B/A'}) = \Omega^1_{\pi_0(B)/\pi_0(A')} = 0$, that $L_{B/A'}$ is of Tor-amplitude $\le 1$).
This follows from the exact triangle $L_{A'/A} \otimes_{A'} B \to L_{B/A} \to L_{B/A'}$.
\end{proof}


\section{Virtual Cartier divisors}
\label{sec:divisors}

\ssec{Virtual Cartier divisors}

\sssec{}
\label{sssec:virtual divisor}
Let $X$ be a derived scheme.
A \emph{virtual effective Cartier divisor} on $X$ is a derived scheme $D$ together with a quasi-smooth closed immersion $i_D : D \hook X$ of virtual codimension $1$.
Thus locally, a virtual effective Cartier divisor on $\Spec(A)$ is of the form $\Spec(A\modmod(f))$ for some $f \in A$.
We will omit the adjective ``effective'' since we do not treat non-effective divisors.

\begin{rem}
Note that for any derived scheme $X$, the collection of virtual Cartier divisors over $X$ forms an \inftyGrpd $\VDiv(X)$.
It can be described as the subgroupoid of $(\DSch_{/X})^\simeq$ whose objects are virtual Cartier divisors $D \hook X$.
Moreover, since the condition of being a virtual Cartier divisor is stable under (derived) base change and fpqc-local, the assignment $X \mapsto \VDiv(X)$ determines a subsheaf of the fpqc sheaf $X \mapsto (\DSch_{/X})^\simeq$.
\end{rem}

\begin{exam}
Suppose that $X$ is classical.
Then any classical (effective) Cartier divisor on $X$ is a virtual Cartier divisor.
\end{exam}

\sssec{}

We record the following computation that will be useful later:

\begin{lem} \label{lem:paths g=0 in A//f}
Let $A$ be a \scr and $f\in A$ an element.
For any point $g \in A$, there is a canonical isomorphism of spaces
  \begin{equation*}
    \Maps_{A\modmod(f)}(g, 0) \simeq \Fib_g(A\xrightarrow{f}A),
  \end{equation*}
between the space of paths $g \simeq 0$ in the underlying space of $A\modmod(f)$, and the space of pairs $(a, \alpha)$, where $a\in A$ is a point and $\alpha : f a \simeq g$ is a path in $A$.
\end{lem}

\begin{proof}
There is a fibre sequence
  \begin{equation*}
    A \xrightarrow{f} A \to A\modmod(f)
  \end{equation*}
of underlying spaces.
\end{proof}

\ssec{Generalized Cartier divisors}

\sssec{}
Let $X$ be a derived scheme.
A \emph{generalized (effective) Cartier divisor} over $X$ is a pair $(\sL, s)$, where $\sL$ is a locally free $\sO_X$-module of rank one and $s$ is an $\sO_X$-module morphism $\sL \to \sO_X$ (see e.g.\ \cite[Def.~10.3.2]{OlssonBook}).

\begin{rem}
According to \cite{IllusieGarden}, generalized Cartier divisors were first introduced by Deligne in a 1988 letter to Illusie, under the name \emph{``divisors''} (where the quotation marks are part of the terminology).
\end{rem}

\sssec{}\label{sssec:generalized divisors}
Any generalized Cartier divisor $(\sL,s)$ on $X$ gives rise to a virtual Cartier divisor as follows.
Denote by $L$ the line bundle $\bV_X(\sL) = \Spec_X(\Sym_{\sO_X}(\sL))$, and let $D$ denote the derived fibred product
  \begin{equation*}
    \begin{tikzcd}
      D \ar{r}{i}\ar{d}
        & X \ar{d}{s}
      \\
      X \ar{r}{0}
        & L
    \end{tikzcd}
  \end{equation*}
where $0$ is the zero section; in other words, $D$ is the derived zero-locus of the section~$s$.
Then $i : D \to X$ is a virtual Cartier divisor with conormal sheaf $\sN_{D/X} = \sL|_D$.

\sssec{}\label{sssec:GDiv}

Let $\GDiv$ denote the derived stack classifying generalized Cartier divisors over~$S$.
To be precise, it is the fpqc sheaf of spaces
  \begin{equation*}
    \GDiv : (\DSch^\aff)^\op \to \Spc,
    \qquad
    S \mapsto (\uPic(S)_{/\sO_S})^\simeq.
  \end{equation*}
Here $\uPic(S)_{/\sO_S}$ is the \inftyCat of pairs $(\sL, s)$ with $\sL \in \uPic(S)$ a locally free sheaf of rank one and $s$ a morphism $\sL \to \sO_S$ (not necessarily invertible), and $(-)^\simeq$ denotes the operation of discarding the non-invertible morphisms.

\sssec{}

Generalized Cartier divisors are classified by the stack $[\A^1/\bG_m]$: the quotient of the affine line $\A^1$ by the canonical $\bG_m$-action by scaling.
This was first observed by L. Lafforgue in 2000 \cite{IllusieGarden}, see also \cite[Prop.~10.3.7]{OlssonBook}.
We next show that this remains valid in derived algebraic geometry, and moreover, that generalized Cartier divisors are the same thing as virtual Cartier divisors.

\begin{prop}\label{prop:VDiv=GDiv}
There are canonical isomorphisms of derived stacks
  \begin{equation*}
    \VDiv \isoto \GDiv \isoto [\A^1/\bG_m].
  \end{equation*}
\end{prop}

\sssec{}

Before proving \propref{prop:VDiv=GDiv} we need to make a brief but slightly technical digression.
Let $D \hook X$ be a virtual Cartier divisor.
Suppose given a quasi-coherent $\sO_X$-algebra $\sA$ and a morphism $\varphi : \sO_D \to \sA$ of $\sO_X$-\emph{modules}.
We would like to show that, if $\varphi$ fits in a commutative triangle
  \begin{equation*}
    \begin{tikzcd}
      \sO_X \ar{r}\ar[swap]{rd}{\eta}
        & \sO_D \ar{d}{\varphi}
      \\
        & \sA,
    \end{tikzcd}
  \end{equation*}
where $\eta$ is the unit, then $\varphi$ lifts to a morphism of $\sO_X$-\emph{algebras} in an essentially unique manner.

More precisely, let $\QcohAlg(X)$ denote the \inftyCat of quasi-coherent $\sO_X$-algebras, and $\Qcoh(X)_{\sO_X\backslash-}$ the \inftyCat of quasi-coherent $\sO_X$-modules $\sF$ equipped with a morphism $\sO_X \to \sF$.
Then we have:

\begin{lem}\label{lem:producing algebra maps A//f -> B}
Let $X$ be a derived scheme and $D \hook X$ a virtual Cartier divisor.
For any quasi-coherent $\sO_X$-algebra $\sA$, the canonical map
  \begin{equation*}
    \Maps_{\QcohAlg(X)}(\sO_D, \sA) \to \Maps_{\Qcoh(X)_{\sO_X\backslash-}}(\sO_D, \sA)
  \end{equation*}
is invertible.
\end{lem}

\begin{proof}
The assignments $X \mapsto \QcohAlg(X)$ and $X \mapsto \Qcoh(X)_{\sO_X\backslash-}$ both form sheaves of \inftyCats, so the question is local on $X$.
Therefore we may assume that $X = \Spec(A)$, for a \scr $A$, and $D = \Spec(A\modmod(f))$ for some $f \in A$.
Note that for any $B \in \SCRing_A$, the space of $A$-algebra homomorphisms $A\modmod(f) \to B$ is the space of commutative squares in $\SCRing_A$
  \begin{equation*}
    \begin{tikzcd}
      A[T] \ar{r}{T\mapsto 0}\ar[swap]{d}{T\mapsto f}
        & A \ar{d}
      \\
      A \ar{r}
        & B
    \end{tikzcd}
  \end{equation*}
as in the proof of \lemref{lem:map-from-A//f}.
This is equivalently the space of commutative triangles in $\Mod_A$
  \begin{equation*}
    \begin{tikzcd}
      A \ar{r}{f}\ar[swap]{rd}{0}
        & A \ar{d}
      \\
        & B
    \end{tikzcd}
  \end{equation*}
which is equivalently the space of $A$-module homomorphisms $A\modmod(f) \to B$ extending $A \to B$, in view of the exact triangle $A \xrightarrow{f} A \to A\modmod(f)$.
\end{proof}

\sssec{Proof of \propref{prop:VDiv=GDiv}}
We first prove that there is an isomorphism $\GDiv \simeq [\A^1/\bG_m]$.
Note that we have a forgetful map $\GDiv\to \uPic^\simeq = B\bG_m$ taking a generalized divisor $(\sL,s)$ to the sheaf $\sL$.
We now consider the derived base change $U\to \GDiv$ of the canonical section $s :\Spec \bZ\to B\bG_m$ that takes a scheme $T$ to $\sO_T$.
The derived fibred product $U$ can be described as the sheaf $S\mapsto (\sL,s,\varphi)$ where $\varphi : \sO_S\to \sL$ is an isomorphism.
But $U\simeq \bA^1$ via the identification $(\sL,s,\varphi)\mapsto \varphi\circ s$.
The resulting map $\bA^1\to \GDiv$ takes a section $f\in \Gamma(S,\sO_S)$ to the generalized divisor $(\sO_S,f)$ and is a $\bG_m$-torsor exhibiting $\GDiv$ as the stack quotient $[\bA^1/\bG_m]$.

We now construct the isomorphism $\VDiv \simeq \GDiv$.
Given an affine derived scheme $S$ and a virtual Cartier divisor $i : D \hook S$, consider the induced morphism $i^\sharp : \sO_S \to i_*(\sO_D)$ and form its homotopy fibre $\sF$, which is equipped with a canonical map $s : \sF \to \sO_S$.
Note that $\sF$ is locally free of rank one, because Zariski-locally on $S$, the exact triangle $\sF \to \sO_S \to i_*(\sO_D)$ takes the form $\sO_S \xrightarrow{f} \sO_S \to \sO_S\modmod(f)$ for some $f \in \Gamma(S, \sO_S)$.
The assignment $(i : D \hook S) \mapsto (\sF, s)$ is clearly functorial and can be regarded as a map of sheaves $p : \VDiv \to \GDiv$.
Observe that the construction of \sssecref{sssec:generalized divisors} shows that this map is surjective on $\pi_0$.

Let $S$ be an affine derived scheme and let $D \hook S$ and $D' \hook S$ be virtual Cartier divisors.
To show that $p$ is invertible it will now suffice to show that the induced map
  \begin{equation*}
    p_{D,D'}: \Maps_{\VDiv(S)}(D, D') \to \Maps_{\GDiv(S)}((\sL, s), (\sL', s'))
  \end{equation*}
is invertible, where $(\sL, s)$ and $(\sL', s')$ are the associated generalized Cartier divisors.
The source of $p_{D,D'}$ is the space of $\sO_S$-algebra isomorphisms $\sO_D \simeq \sO_{D'}$.
The target is the space of commutative squares of $\sO_S$-modules
  \begin{equation*}
    \begin{tikzcd}
      \sL \ar{r}{s}\ar{d}{u}
        & \sO_S \ar[equals]{d}
      \\
      \sL' \ar{r}{s'}
        & \sO_S
    \end{tikzcd}
  \end{equation*}
where $u$ is an isomorphism, which is tautologically equivalent to the space of $\sO_S$-module isomorphisms $\sO_D \simeq \sO_{D'}$ compatible with the maps $\sO_S \to \sO_D$ and $\sO_S \to \sO_{D'}$ (via a specified commutative triangle).
The fact that $p_{D,D'}$ is invertible is the content of \lemref{lem:producing algebra maps A//f -> B}.


\section{Blow-ups}
\label{sec:blow-up}


\ssec{Construction and main properties}
\label{ssec:construction}

\sssec{}
\label{sssec:blow-up}

Let $i : Z \hook X$ be a quasi-smooth closed immersion of derived schemes.
For any derived scheme $S$ and morphism $f : S \to X$, a \emph{virtual Cartier divisor on $S$ lying over} $(X,Z)$ is the datum of a commutative square
  \begin{equation} \label{eq:virtual divisor lying over Z}
    \begin{tikzcd}
      D \ar{d}{g}\ar[hookrightarrow]{r}{i_D}
        & S \ar{d}{f}
      \\
      Z \ar[hookrightarrow]{r}{i}
        & X
    \end{tikzcd}
  \end{equation}
satisfying the following conditions:

\begin{enumerate}
\item 
The morphism $i_D : D \to S$ exhibits $D$ as a virtual Cartier divisor on $S$.

\item
The underlying square of classical schemes is cartesian.

\item
The canonical morphism
  \begin{equation} \label{eq:comparison map of conormal sheaves}
    g^*\sN_{Z/X} \to \sN_{D/S}
  \end{equation}
is surjective (on $\pi_0$).
\end{enumerate}

\begin{exam}
Suppose that $X$, $Z$, and $S$ are classical schemes.
If the classical schematic fibre $f^{-1}(Z)$ is a classical Cartier divisor on $S$, then $f^{-1}(Z)$ also defines a virtual Cartier divisor lying over $(X,Z)$.
Condition (c) follows from \cite[Cor.~25.3.6.4]{SAG-20180204}, also see \remref{rem:S_Z}\itemref{rem:S_Z/S_Z-props-wo-cotangent-complex}.
\end{exam}

\begin{rems}\label{rem:S_Z}
Let $S_Z$ denote the derived fibred product $S \fibprod^\bR_X Z$.
\begin{remlist}
\item
If the square \eqref{eq:virtual divisor lying over Z} is homotopy cartesian, that is, if $D\simeq S_Z$, then the morphism \eqref{eq:comparison map of conormal sheaves} is an isomorphism.
If $i$ is of virtual codimension $n > 1$ then this is never the case.

\item\label{rem:S_Z/S_Z-props}
We can think of a virtual Cartier divisor on $S$ lying over $(X,Z)$ equivalently as a derived scheme $D$ over $S_Z$ such that (a) the induced morphism $i_D : D \to S_Z \hook S$ is a virtual Cartier divisor; (b) the morphism $D \to S_Z$ induces an isomorphism $D_\cl \simeq (S_Z)_\cl$ on underlying classical schemes; and (c) the canonical morphism $h^*\sN_{S_Z/S} \to \sN_{D/S}$ is surjective, where $h : D \to S_Z$.
This latter condition is equivalent to the relative cotangent complex $\sL_{D/S_Z}$ being $2$-connective ($\pi_{i\le 1} = 0$).

\item\label{rem:S_Z/S_Z-props-wo-cotangent-complex}
The closed immersion $h : D \to S_Z$ induces a map $\phi : \sO_{S_Z}\to h_*\sO_D$ which always is surjective on $\pi_0$.
Condition (b) is equivalent to $\pi_0(\phi)$ being an isomorphism.
Since $\pi_1(L_{D/S_Z})=\pi_0(\Fib(h^*\phi))$ by the connectivity properties of the Hurewicz map $\epsilon_\phi$ \cite[Prop.~25.3.6.1]{SAG-20180204}, we see that condition (c) is equivalent to the surjectivity of $\pi_1(h^*\phi)$.
By Nakayama, this is equivalent to the surjectivity of $\pi_1(\phi)$ when $|D|=|S_Z|$, i.e., under condition (b).
We conclude that (b)+(c) is equivalent to: $\pi_0(\phi)$ is an isomorphism and $\pi_1(\phi)$ is surjective.
\end{remlist}
\end{rems}

\sssec{}
\label{sssec:Bl_Z/X}

The collection of virtual Cartier divisors on $S$ lying over $(X,Z)$ forms a space which we denote $\Bl_{Z} X(S \to X)$.
Moreover, the construction is functorial and defines a presheaf of spaces
  \begin{equation*}
    \Bl_{Z} X : (\DSch_{/X})^\op \to \Spc
  \end{equation*}
on the site of derived schemes over $X$.
Indeed consider the presheaf $\sF : (S\to X)\mapsto(\DSch_{/S_Z})^\simeq$, which sends $S \to X$ to the space obtained by discarding non-invertible morphisms in the \inftyCat $\DSch_{/S_Z}$, and note that $\Bl_{Z} X$ defines a sub-presheaf since the conditions (a), (b) and (c) are stable under derived base change.
Since these conditions are also étale-local, and $\sF$ satisfies (hyper)descent for the étale topology, the presheaf $\Bl_{Z} X$ is also an étale (hyper)sheaf.
In particular, it defines a derived stack $\Bl_{Z} X$ over $X$, which we call the \emph{blow-up}; we denote by $\pi_{Z/X}$ the structural morphism $\Bl_{Z} X \to X$.
The main properties of the construction $\Bl_{Z} X$ are summed up below:

\begin{thm}\leavevmode\label{thm:all}
Let $i : Z\hook X$ be a quasi-smooth closed immersion of derived schemes.
\begin{thmlist}
\item\label{thm:all/representable}
The derived stack $\Bl_{Z} X$ is (representable by) a derived scheme.

\item\label{thm:all/stable under base change}
The construction $\Bl_{Z} X \to X$ commutes with arbitrary derived base change.
That is, $(\Bl_{Z} X) \times^{\bR}_X X'=\Bl_{Z\times^{\bR}_X X'} X'$ for every morphism $X'\to X$ of derived schemes.

\item\label{thm:all/covariance}
The construction $\Bl_{Z} X \to X$ has covariant functoriality in $X$ along quasi-smooth closed immersions.
That is, for any quasi-smooth closed immersion $X \to Y$, there is a canonical quasi-smooth closed immersion
  \begin{equation*}
    \Bl_{Z} X \to \Bl_{Z} Y
  \end{equation*}
of derived schemes over $Y$.

\item\label{thm:all/exceptional}
There is a canonical closed immersion $\P_Z(\sN_{Z/X}) \hook \Bl_{Z} X$ which exhibits the projectivized normal bundle as the universal virtual Cartier divisor lying over $(X,Z)$.

\item\label{thm:all/proper quasi-smooth}
The structural morphism $\pi_{Z/X} : \Bl_{Z} X \to X$ is proper and quasi-smooth, and induces an isomorphism $\Bl_{Z} X \setminus \P_Z(\sN_{Z/X}) \isoto X\setminus Z$.

\item\label{thm:all/classical}
Suppose that $X$ and $Z$ are classical schemes.
Then the derived scheme $\Bl_{Z} X$ is classical, and coincides with the classical blow-up $\Bl^\cl_{Z} X$.

\item\label{thm:all/classical-description}
In general, $\left(\Bl_{Z} X\right)_\cl=\P^\cl_{X_\cl}\bigl(\pi_0(\sI)\bigr)$ where $\sI=\Fib(\sO_X\to i_*\sO_Z)$ is the homotopy fibre and $\P^\cl(-)=\Proj(\Sym(-))$ is the homogeneous spectrum of the underived symmetric algebra.

\item\label{thm:all/divisor}
If $i$ is a virtual Cartier divisor, then the morphism $\pi_{Z/X} : \Bl_{Z} X \to X$ is invertible.

\item\label{thm:all/id}
If $i = \id_X$, then the blow-up $\Bl_{X} X$ is empty.
\end{thmlist}
\end{thm}

The proof will be delayed a few pages (see \ssecref{ssec:proof}).

\sssec{}

\thmref{thm:all} admits the following immediate generalization.
Given a closed immersion $i : Z \hook X$ of derived stacks, we say that $i$ is quasi-smooth if, fpqc-locally on $X$, it is a quasi-smooth closed immersion of derived schemes.
Then the discussion above applies \emph{mutatis mutandis} to define a derived stack $\Bl_{Z} X$, classifying virtual Cartier divisors lying over $(X,Z)$.
Moreover, if $X$ is an ($n$-geometric) derived Deligne--Mumford (resp.\ Artin) stack in the sense of \cite{HAG2}, then the same is true of $\Bl_{Z} X$.

\sssec{}

The following universal property for the classical blow-up $\Bl^\cl_{Z} X$ follows from part \itemref{thm:all/classical} and the definition of $\Bl_{Z} X$:

\begin{cor} \label{cor:universal property in classical case}
Let $i : Z \hook X$ be a regular closed immersion between classical schemes.
For any classical scheme $S$ over $X$, the set of $X$-morphisms $S \to \Bl^\cl_{Z} X$ is in bijection with the set of virtual Cartier divisors on $S$ lying over $(X,Z)$.
\end{cor}

\begin{rem}
In the situation of the Corollary, assertion \itemref{thm:all/classical} in the Theorem implies that for any classical scheme $S \to X$, the space $\Bl_{Z} X(S \to X)$ of virtual Cartier divisors on $S$ lying over $(X,Z)$ is \emph{discrete} (i.e., can indeed be identified with a set).
\end{rem}

\sssec{}
\label{sssec:KST}

Suppose that $X = \Spec(A)$ is a classical noetherian affine scheme, $Z = \Spec(A/I)$ is a closed subscheme, and $f_1,\ldots,f_n$ are generators of the ideal $I$.
In this situation the authors of \cite{KST} consider the derived scheme
  \begin{equation*}
    X \fibprod^\bR_{\A^n} \Bl^\cl_{\{0\}} \A^n,
  \end{equation*}
where $f : X \to \A^n$ corresponds to the elements $f_i$.
Parts \itemref{thm:all/stable under base change} and \itemref{thm:all/classical} of the Theorem above show that this derived scheme is nothing else than the blow-up of $\Spec(A)$ in the quasi-smooth derived subscheme $\tilde{Z} = \Spec(A\modmod(f_1,\ldots,f_n))$.
In particular it follows that the construction of \emph{op.\ cit.}\ is intrinsic not to the elements $f_i$, but only to their derived zero-locus $\tilde{Z}$.

\sssec{}
Let $X$ be a classical scheme and $Z$ a finitely presented closed subscheme.
Although the ideal $\sI$ defining $Z$ is not always globally generated, we can often find a vector bundle $\sE$ and a surjection $\varphi : \sE\to \sI$.
This is for example the case when $X$ is quasi-projective or, more generally, has the resolution property.
The data $(\sE,\varphi)$ endows $Z$ with a quasi-smooth derived structure $W$ as follows.
The induced morphism $\sE\to \sI\to \sO_X$ corresponds to a section $s$ of the vector bundle $E=\bV_X(\sE) = \Spec_X(\Sym_{\sO_X}(\sE))$.
Consider the derived fibred product
  \begin{equation*}
    \begin{tikzcd}
      W \ar{r}{i}\ar{d} & X \ar{d}{s} \\
      X \ar{r}{0} & E.
    \end{tikzcd}
  \end{equation*}
Then $W\to X$ is a quasi-smooth closed immersion and $W_\cl=Z$.
Note that the classical blow-up $\Bl^\cl_{Z} X$ need not equal the underlying classical scheme of the derived blow-up $\Bl_{W} X$ but is always the schematic closure of $X\smallsetminus Z$ in $\Bl_{X} E \times_E X = (\Bl_{W} X)_\cl$.

\sssec{}

An immediate application of \thmref{thm:all} is the existence of a \emph{deformation to the normal bundle}, for any quasi-smooth closed immersion.
This was first constructed by Verdier \cite{VerdierRR} for regular closed immersions of classical schemes.
The construction gives a deformation of a quasi-smooth closed immersion $Z\hook X$ to the zero-section $Z\hook \mrm{N}_{Z/X}$.

\begin{thm}\label{thm:deformation}
Let $i : Z \hook X$ be a quasi-smooth closed immersion of derived stacks with normal bundle $\mrm{N}_{Z/X} = \Spec_Z(\Sym_{\sO_Z}(\sN_{Z/X}))$.
Then there exists a canonical factorization of $i\times \id$, the \emph{deformation to the normal bundle}:
\[
i\times\id : Z\times \A^1 \xhookrightarrow{\;j\;} \mrm{D}_{Z/X} \xrightarrow{\;\pi\;} X\times \A^1
\]
satisfying the following properties
\begin{thmlist}
\item
The factorization is stable under arbitrary derived base change along $X$.
\item $j$ is a quasi-smooth closed immersion.
\item $\pi$ is quasi-smooth and quasi-projective.
\item Restricting to $\bG_m=\A^1\setminus \{0\}$ we obtain
\[
i\times \id : Z\times \bG_m \xhookrightarrow{\;j_\mrm{gen}\;} X\times \bG_m \xrightarrow{\;\pi_\mrm{gen}\;} X\times \bG_m
\]
where $j_\mrm{gen}=i\times \id$ and $\pi_\mrm{gen}=\id$.
\item\label{thm:deformation/special}
Restricting to $\{0\}$ we obtain
\[
i : Z \xhookrightarrow{\;j_0\;} \mrm{N}_{Z/X} \xrightarrow{\;\pi_0\;} X
\]
where $j_0$ is the zero-section and $\pi_0$ the composition of the projection $\mrm{N}_{Z/X}\to Z$ and $i : Z\hook X$.
\end{thmlist}
In particular, we have a virtual Cartier divisor $\mrm{N}_{Z/X}\hook \mrm{D}_{Z/X}$.
\end{thm}
\begin{proof}
We begin by considering the following \emph{projective} variant:
\[
\overline{\pi} : \overline{\mrm{D}}_{Z/X}=\Bl_{Z\times\{0\}} (X\times \A^1) \to X\times \A^1.
\]
We have a quasi-smooth closed immersion
\[
\overline{j} : Z\times\A^1=\Bl_{Z\times\{0\}} (Z\times\A^1) \hook
\Bl_{Z\times\{0\}} (X\times\A^1) =\overline{\mrm{D}}_{Z/X}
\]
by \thmref{thm:all} \itemref{thm:all/covariance} and \itemref{thm:all/divisor}.
Since
\[
Z\times\{0\}=(Z\times\A^1)\times_{X\times\A^1} (X\times \{0\})
\]
we have that $\Bl_{Z\times\{0\}} (Z\times\A^1)$ and $\Bl_{Z\times\{0\}} (X\times\{0\})$ are disjoint inside $\overline{\mrm{D}}_{Z/X}$.
We thus define
\[
\mrm{D}_{Z/X}=\overline{\mrm{D}}_{Z/X}\setminus \Bl_{Z\times\{0\}} (X\times\{0\})
\]
and let $j$ and $\pi$ be the restrictions of $\overline{j}$ and $\overline\pi$.
It remains to prove \itemref{thm:deformation/special}.
The special fibre of $\mrm{D}_{Z/X}$ is the exceptional divisor of the blow-up of $X\times\A^1$ minus the exceptional divisor of the blow-up of $X\times \{0\}$, that is:
\[
\mrm{D}_{Z/X}\times_{\A^1} \{0\} = \P_Z(\sN_{Z\times\{0\}/X\times \A^1}) \setminus \P_Z(\sN_{Z/X})
\]
As $\sN_{Z\times\{0\}/X\times\A^1}$ splits as the direct sum $\sN_{Z/X} \oplus \sO_{Z}$, this is $\mrm{N}_{Z/X}$.
\end{proof}

Combining \thmref{thm:deformation} with Fulton's methods in \cite{Fulton}, one can construct virtual Gysin maps along quasi-smooth closed immersions in Chow theory.
This gives in particular a simple new construction of virtual fundamental classes, the details of which will be given elsewhere.


\ssec{A special case}
\label{ssec:special case}

In this subsection we will study the special case of the quasi-smooth closed immersion $\{0\} \hook \A^n$, $n\ge 1$.
We will identify $\A^n$ with $\Spec(\bZ[T_1,\ldots,T_n])$ and write $Y := \Bl_{\{0\}} \A^n$ for convenience.

\sssec{}
\label{sssec:main computation}

For each $1\le k\le n$ let $A_k = \bZ[T_1/T_k,\ldots,T_n/T_k,T_k]$.
The commutative squares
  \begin{equation*}
    \begin{tikzcd}
      \Spec(A_k/(T_k)) \ar[hookrightarrow]{r}\ar{d}
        & \Spec(A_k)\ar{d}
      \\
      \Spec(\bZ[T_1,\ldots,T_n]/(T_1,\ldots,T_n))\ar[hookrightarrow]{r}
        & \Spec(\bZ[T_1,\ldots,T_n])
    \end{tikzcd}
  \end{equation*}
define virtual Cartier divisors lying over $(\A^n,\{0\})$, which are classified by canonical morphisms $\Spec(A_k) \to Y$.

\begin{lem}\label{lem:special}\leavevmode
\begin{thmlist}
\item\label{lem:special/charts-mono}
For each $1\le k\le n$, the morphism $\Spec(A_k) \to Y$ is a monomorphism.
That is, its homotopy fibres are empty or contractible.

\item\label{lem:special/charts-epi}
The induced morphism $\coprod_k \Spec(A_k) \to Y$ is an effective epimorphism.
In other words, the family $(Y_k \hook Y)_k$ defines a Zariski atlas for the derived stack $Y$.

\item\label{lem:special/classical}
The derived stack $Y$ is a classical scheme.
Moreover, it is isomorphic to the classical blow-up $\Bl^\cl_{\{0\}} \A^n$.
\end{thmlist}
\end{lem}

\sssec{Proof of \lemref{lem:special}\itemref{lem:special/charts-mono}}

Let $S = \Spec(R)$ be an affine derived scheme and $f : S \to \A^n$ a morphism corresponding to points $f_1,\ldots,f_n \in R$.
It suffices to show that the induced map of spaces
  \begin{equation*}
    \theta : \Maps_{{/\A^n}}(S, \Spec(A_k)) \to Y_k(S \to \A^n)
  \end{equation*}
is a monomorphism.
Set $A := \Maps_{{/\A^n}}(S, \Spec(A_k))$ and $B := Y_k(S \to \A^n)$.
We will view $A$ and $B$ as \inftyGrpds and show that $\theta$ is fully faithful.

First observe that the source $A$ can be described as
  \begin{equation*}
    A = \Maps_{\SCRing_{\bZ[T_1,\ldots,T_n]}}(A_k, R) \simeq \prod_{r\ne k} \Fib_{f_r}(R \xrightarrow{f_k} R),
  \end{equation*}
since $A_k = \bZ[T_1/T_k,\ldots,T_n/T_k,T_k]$ is free, as a simplicial commutative $\bZ[T_1,\ldots,T_n]$-algebra, on generators $X_1,\ldots,\hat{X_k},\ldots,X_n$ with relations $T_k X_r-T_r$, $r \ne k$ (because the sequence $(T_k X_r-T_r)_{r\ne k}$ is regular, see the proof of \cite[Exp.~VII, Prop.~1.8(ii)]{SGA6}).
Thus its objects can be identified with tuples $(a_r,\alpha_r)_{r\ne k}$, where $a_r \in R$ are points and $\alpha_r : f_k a_r \simeq f_r$ are paths in (the underlying space of) $R$.

The target $B$ is by definition a sub-\inftyGrpd of $(\DSch^{\aff}_{/S\fibprod_{\A^n}^\bR \{0\}})^\simeq$, that is, $(\SCRing_{R\modmod(f_1,\ldots,f_n)})^\simeq$, where the notation $(-)^\simeq$ means that we take the sub-\inftyGrpd of the \inftyCat consisting of only invertible morphisms.
Thus its objects are morphisms $R\modmod(f_1,\ldots,f_n) \to R'$, where $R'$ is a \scr.
These are equivalently $R$-algebras $R'$ equipped with paths $f_j \simeq 0$ in $R'$, $1\le j\le n$.

The map $\theta : A \to B$ sends an object $a = (a_r,\alpha_r)_{r\ne k} \in A$ to the object $\theta(a) \in B$ given by the $R$-algebra $R\modmod(f_k)$ together with certain paths $\theta(a)_j : f_j \simeq 0$ in $R\modmod(f_k)$, $1\le j\le n$.
The path $\theta(a)_k$ is the ``tautological'' path and the other paths $\theta(a)_r$, $r\ne k$, are induced by composing $\alpha_r$ with $\theta(a)_k$.
Let $a'$ be another object of $A$ and $\theta(a')$ its image in $B$.
The space of paths $\theta(a) \simeq \theta(a')$ can be described as follows:
  \begin{align*}
    \Maps_{B}(\theta(a), \theta(a'))
      &= \Maps_{\SCRing_{R\modmod(f_1,\ldots,f_n)}}(R\modmod(f_k), R\modmod(f_k))\\
      \begin{split}
         &= \Fib_{\theta(a')}\Bigl(\Maps_{\SCRing_{R}}(R\modmod(f_k), R\modmod(f_k))\\
          &\quad\qquad\xrightarrow{\theta(a)^*} \Maps_{\SCRing_{R}}(R\modmod(f_1,\ldots,f_n), R\modmod(f_k))\Bigr)
        \end{split}\\
      &= \Fib_{(\theta(a')_j)_j}\Bigl(\Maps_{R\modmod(f_k)}(f_k, 0) \to \prod_{j=1}^n \Maps_{R\modmod(f_k)}(f_j, 0)\Bigr)\\
      &= \prod_{r\ne k} \Maps_{R\modmod(f_k)}(\theta(a)_r, \theta(a')_r).
  \end{align*}
Under the above identifications the map $\Maps_A(a,a') \to \Maps_B(\theta(a),\theta(a'))$ is identified with the canonical map
  \begin{equation*}
    \prod_{r\ne k} \Maps_{\Fib_{f_r}(R \xrightarrow{f_k} R)}((a_r,\alpha_r),(a'_r,\alpha'_r)) \to \prod_{r\ne k} \Maps_{R\modmod(f_k)}(\theta(a)_r, \theta(a')_r)
  \end{equation*}
which is invertible by \lemref{lem:paths g=0 in A//f}.

\sssec{Proof of \lemref{lem:special}\itemref{lem:special/charts-epi}}

Let $S$ be a derived scheme, $f : S \to \A^n$ a morphism, and $D$ a virtual Cartier divisor on $S$ lying over $(\A^n,\{0\})$.
The claim is that Zariski-locally on $S$, $D$ fits into a cartesian square
  \begin{equation} \label{eq:Y_k}
    \begin{tikzcd}
      D \ar[hookrightarrow]{r}\ar{d}
        & S \ar{d}
      \\
      \Spec(A_k/(T_k)) \ar[hookrightarrow]{r}\ar{d}
        & \Spec(A_k) \ar{d}
      \\
      \{0\} \ar[hookrightarrow]{r}
        & \A^n
    \end{tikzcd}
  \end{equation}
for some $k$.
We can assume that $S$ is affine, say $S = \Spec(R)$ for some $R \in \SCRing$, so that $f$ corresponds to points $f_1,\ldots,f_n \in R$.
The fact that $D$ lies over $(\A^n,\{0\})$ then implies that locally, the conormal sheaf $\sN_{D/S}$ has a basis given by $df_k$ for some $k$.
It follows that the induced morphism $D \to \Spec(R\modmod(f_1,\ldots,f_n))\to \Spec(R\modmod(f_k))$ is invertible, arguing as in the end of the proof of \propref{prop:regular immersions and cotangent}.

\sssec{Proof of \lemref{lem:special}\itemref{lem:special/classical}}

Assertions~\itemref{lem:special/charts-mono} and~\itemref{lem:special/charts-epi} provide a Zariski cover for $Y$ by the standard affine cover of $\Bl^\cl_{\{0\}} \A^n$, so the claim follows.

\begin{rem} \label{rem:classical universal property}
Once one knows that $Y$ is classical scheme, it is easy to check directly that $Y$ satisfies the classical universal property of $\Bl^\cl_{\{0\}} \A^n$.
That is, suppose $S$ is a classical scheme and $f : S \to \A^n$ is a morphism.
If the classical schematic fibre $f^{-1}(\{0\})$ is a classical Cartier divisor on $S$, then it lies over $(\A^n,\{0\})$ as a virtual Cartier divisor, and is moreover the \emph{unique} such; in particular, there exists a unique morphism $S \to Y$ over $\A^n$.
\end{rem}


\ssec{Proof of main theorem (\sssecref{sssec:Bl_Z/X})}
\label{ssec:proof}

\sssec{Proof of \ref{thm:all/stable under base change}}
\label{sssec:proof of stable under base change}

Let $i : Z \hook X$ be a quasi-smooth closed immersion of derived schemes, and $i' : Z' \hook X'$ its derived base change along a morphism $p : X' \to X$.
Given a derived scheme $S'$ over $X'$, any virtual Cartier divisor $D'$ on $S'$ lying over $(X',Z')$ also lies over $(X,Z)$.
In particular there is a canonical morphism $\Bl_{Z'} X' \to \Bl_{Z} X \fibprod^\bR_X X'$.
We prove the following more precise formulation of the statement:

\begin{claim} \label{claim:base change}
The canonical morphism of derived stacks $\Bl_{Z'} X' \to \Bl_{Z} X \fibprod^\bR_X X'$ is invertible.
\end{claim}

\begin{proof}
Use the description mentioned in \remref{rem:S_Z}\itemref{rem:S_Z/S_Z-props}: for any $S' \to X'$, the spaces $\Bl_{Z'} X'(S'\to X')$ and $(\Bl_{Z} X \fibprod^\bR_X X')(S' \to X')$ both define the same subspace of $(\DSch_{/S'\fibprod^\bR_{X'} Z'})^\simeq = (\DSch_{/S'\fibprod^\bR_{X} Z})^\simeq$.
\end{proof}

\sssec{Proof of \ref{thm:all/representable}}

It suffices to show this Zariski-locally on the base $X$, so we can assume that $i : Z \hook X$ is a derived base change of $\{0\} \hook \A^n$.
Derived fibred products of derived schemes are representable, so by \thmref{thm:all}\ref{thm:all/stable under base change} (proven in \sssecref{sssec:proof of stable under base change} above) we can reduce to the special case considered in \lemref{lem:special}.

\sssec{Proof of \ref{thm:all/exceptional}}

Let $D^\univ_{Z/X}$ denote the ``universal virtual Cartier divisor'' lying over $(X,Z)$, classified by the identity morphism $\Bl_{Z} X \to \Bl_{Z} X$.
This is a derived scheme $D^\univ_{Z/X}$ equipped with a canonical morphism $\pi_\univ : D^\univ_{Z/X} \to Z$, a canonical locally free sheaf $\sL^\univ_{Z/X} := \sN_{D^\univ_{Z/X}/\Bl_{Z} X}$ of rank $1$, and a canonical surjection $(\pi_\univ)^*\sN_{Z/X} \to \sL^\univ_{Z/X}$.
This data is classified by a canonical morphism
  \begin{equation*}
    D^\univ_{Z/X} \to \P_Z(\sN_{Z/X})
  \end{equation*}
of derived schemes over $Z$.

\begin{claim}
The morphism $D^\univ_{Z/X} \to \P_Z(\sN_{Z/X})$ is invertible.
In particular, there is a canonical closed immersion $\P_Z(\sN_{Z/X}) \hook \Bl_{Z} X$ which exhibits the projectivized normal bundle as the universal virtual Cartier divisor lying over $(X,Z)$.
\end{claim}

\begin{proof}
The assertion is local and stable under derived base change, so we reduce to the case of $\{0\} \hook \A^n$.
Then we can apply the well-known universal property of the classical blow-up (\remref{rem:classical universal property}): since the classical fibred product $\Bl_{\{0\}} \A^n \fibprod_{\A^n} \{0\}$ is the classical effective Cartier divisor $\P_{\{0\}}(\sN_{\{0\}/\A^n})$, the conclusion is that there is a unique virtual Cartier divisor $\P_{\{0\}}(\sN_{\{0\}/\A^n}) \hook \Bl_{\{0\}} \A^n$ lying over $(\A^n,\{0\})$, classified by the identity of $\Bl_{\{0\}} \A^n$.
\end{proof}

\sssec{Proof of \ref{thm:all/proper quasi-smooth}}

The properties in question are Zariski-local on the target and stable under arbitrary derived base change, so we again reduce to the case of $\{0\} \hook \A^n$.
Then these are well-known properties of the classical blow-up.
In fact, the projection $\Bl_{\{0\}} \A^n \to \A^n$ factors through a regular closed immersion $\Bl_{\{0\}} \A^n \hook \P^{n-1}_{\A^n}$, see \cite[Exp.~VII, Prop.~1.8(ii)]{SGA6}.

\sssec{Proof of \ref{thm:all/classical}}

Suppose that $X$ and $Z$ are classical schemes.
To show that $\Bl_{Z} X$ is classical, we can assume that $X = \Spec(R)$ and $Z = \Spec(R\modmod(f_1,\ldots,f_n)) = \Spec(R/(f_1,\ldots,f_n))$, where $R$ is a commutative ring and $(f_1,\ldots,f_n)$ is a regular sequence.
Then by \thmref{thm:all}\ref{thm:all/stable under base change} and \lemref{lem:special}, the derived scheme $\Bl_{Z} X$ admits a Zariski cover by the schemes
  \begin{equation*}
    \Spec(R) \fibprod^\bR_{\A^n} \Spec(A_k) = \Spec(R[X_1,\ldots,\hat{X_k},\ldots,X_n]\modmod(X_rf_i - f_r)_{r\ne j}),
  \end{equation*}
where $A_k = \bZ[T_1/T_k,\ldots,T_n/T_k,T_k]$ as in \ssecref{ssec:special case}.
Thus the claim follows from the fact that the sequence $(X_rf_i - f_r)_{r\ne j}$ is regular (see the proof of \cite[Exp.~VII, Prop.~1.8(ii)]{SGA6}).

To show that $\Bl_{Z} X$ moreover coincides with $\Bl^\cl_{Z} X$, we can assume $Z \hook X$ is a derived base change of $\{0\} \hook \A^n$ along some morphism $f : X \to \A^n$.
Then we have canonical isomorphisms $\Bl_{Z} X \simeq \Bl_{\{0\}} \A^n \fibprod^\bR_{\A^n} X \simeq \Bl_{\{0\}} \A^n \fibprod_{\A^n} X$ by \thmref{thm:all}\ref{thm:all/stable under base change} and the first part of \ref{thm:all/classical}.
On the other hand the classical blow-up $\Bl^\cl_{Z} X$ is the classical base change $\Bl^\cl_{\{0\}} \A^n \fibprod_{\A^n} X$ by \cite[Exp.~VII, Prop.~1.8(i)]{SGA6}.
Therefore the claim follows from the special case of $\{0\}\hook \A^n$ (\lemref{lem:special}).

\sssec{Proof of \ref{thm:all/classical-description}}
Since the assignment $(i: Z \to X) \mapsto \pi_0\Fib(\sO_X\to i_*\sO_Z)$ 
and the underived $\Sym$ and $\Proj$ commute with underived base change, it is enough to prove that $\Bl^\cl_{\{0\}} \A^n=\Proj(\Sym(I))$ where $I=(T_1,T_2,\dots,T_n) \subset \bZ[T_1,\ldots,T_n]$.
But $I$ is regular so $\Sym(I)$ coincides with the Rees algebra $R(I)$~\cite[Ch.~1, Thm.~1]{Micali}.
The result follows.

\begin{rem}
The blow-up $\Bl_{Z} X$ coincides with $\Proj(\bL\Sym(I))$ if and only if the virtual codimension of $Z\hook X$ is at most $2$.
Indeed, it is enough to consider the situation of $Z=\{0\}$ in $X=\A^n$ so $I=(T_1,T_2,\dots,T_n)$.
An explicit calculation then shows that $\bL\Sym(I)=\Sym(I)[0]$ when $n\geq 2$.
One also calculates that the cotangent complex of the special fibre of $\Proj(\bL\Sym(I))$ is perfect of Tor-amplitude $[0,n-1]$ so $\Proj(\bL\Sym(I))\to X$ is not quasi-smooth when $n\geq 3$, hence cannot equal the quasi-smooth morphism $\Proj(\Sym(I))\to X$.
\end{rem}

\sssec{Proof of \ref{thm:all/divisor}}

Suppose that $i : Z \hook X$ is of virtual codimension $1$.
Then given a derived scheme $S$ over $X$, a virtual Cartier divisor $D \hook S$ lies over $(X, Z)$ if and only if the square \eqref{eq:virtual divisor lying over Z} is homotopy cartesian.
Indeed, if the canonical morphism \eqref{eq:comparison map of conormal sheaves} is surjective then it is an isomorphism and $\sL_{D/S_Z}=0$ so $D\to S_Z$ is an isomorphism.
In other words, $(X, Z)$ is in this case itself the universal virtual Cartier divisor lying over $(X,Z)$.

\sssec{Proof of \ref{thm:all/id}}

Suppose that $i$ is the identity of $X$.
Then given a derived scheme $S$ over $X$, a virtual Cartier divisor $D \hook S$ lies over $(X, X)$ if and only if $S = D = \initial$.
Indeed, the surjectivity of \eqref{eq:comparison map of conormal sheaves} implies that $D=\initial$ and hence that $X=\initial$.

\sssec{Proof of \ref{thm:all/covariance}}

The existence of the canonical morphism $\Bl_{Z} X \to \Bl_{Z} Y$ follows from the observation that any virtual Cartier divisor lying over $(X, Z)$ also lies over $(Y, Z)$.
In other words, this morphism is classified by the commutative square
  \begin{equation*}
    \begin{tikzcd}
      \P_Z(\sN_{Z/X}) \ar{r}\ar{d}
        & \Bl_{Z} X \ar{d}
      \\
      Z \ar{r}
        & Y,
    \end{tikzcd}
  \end{equation*}
viewed as a virtual Cartier divisor lying over $(Z,Y)$.
To show that it is a quasi-smooth closed immersion, we can use \itemref{thm:all/proper quasi-smooth}: over $Y \setminus Z$ it induces the quasi-smooth closed immersion $X \setminus Z \to Y \setminus Z$, and over $Z$ it induces a closed immersion $\P(\sN_{Z/X}) \to \P(\sN_{Z/Y})$ between smooth schemes over $Z$, hence is a quasi-smooth closed immersion.


\section{Simultaneous blow-up in multiple centres}
\label{sec:multiple}

\ssec{}
We begin by reviewing the traditional approach to blowing up several centres on classical schemes.

\sssec{}
Let $X$ be a smooth scheme and let $Z_1$ and $Z_2$ be smooth closed subschemes.
In many situations, one wants to blow-up both $Z_1$ and $Z_2$.
More precisely, one wants a smooth scheme dominating both $\Bl_{Z_1} X$ and $\Bl_{Z_2} X$.
Consider the following three blow-up procedures:
\begin{enumerate}
\item First blow-up $Z_1$, then blow-up the strict transform of $Z_2$.
\item First blow-up $Z_2$, then blow-up the strict transform of $Z_1$.
\item First blow-up $Z_1\cap Z_2$, then blow-up the strict transform of $Z_1\cup Z_2$ which is now the disjoint union of the strict transforms of $Z_1$ and $Z_2$.
\end{enumerate}
For these blow-ups to be smooth, one needs $Z_1\cap Z_2$ to be smooth which can be done using embedded resolution of singularities in characteristic zero.

Let $X_a$, $X_b$ and $X_c$ be the results of these three procedures.
In general, the $X_a$ and $X_b$ are non-isomorphic and do not dominate both $\Bl_{Z_1} X$ and $\Bl_{Z_2} X$ whereas $X_c$ dominates $X_a$, $X_b$, $\Bl_{Z_1} X$ and $\Bl_{Z_2} X$.

\sssec{}
When $Z_1\cap Z_2$ is smooth, then the following are equivalent.
\begin{enumerate}
\item $X_a=X_b$.
\item $Z_1$ and $Z_2$ meets transversely, that is, $\codim(Z_1,X)+\codim(Z_2,X)=\codim(Z_1\cap Z_2,X)$.
\item $Z_1\to X$ and $Z_2\to X$ are Tor-independent, that is, the derived intersection $Z_1\fibprod^\bR_X Z_2$ is a classical scheme.
\end{enumerate}
When these conditions hold, then $X_a=X_b=\Bl_{Z_1} X\fibprod^\bR_X \Bl_{Z_2} X=\Bl_{Z_1} X\times_X \Bl_{Z_2} X$.

\sssec{}
In derived algebraic geometry, it is not a problem if $Z_1\cap Z_2$ is singular or if the intersection is not transversal.
Given a quasi-smooth derived scheme $X$ and two quasi-smooth closed immersions $Z_1\hookrightarrow X$ and $Z_2\hookrightarrow X$ the derived fibred product $\Bl_{Z_1} X\fibprod^\bR_X \Bl_{Z_2} X$ is a quasi-smooth derived scheme that dominates both $\Bl_{Z_1} X$ and $\Bl_{Z_2} X$.

The goal of this section is a construction where the simultaneous blow-up $\Bl_{Z_1} X\fibprod^\bR_X \Bl_{Z_2} X$ of $Z_1$ and $Z_2$ is described as the blow-up of $X$ in $Z_1\amalg Z_2$.

\begin{exam}
If $Z_1\to X$ and $Z_2\to X$ are Tor-independent, then so are $\Bl_{Z_1} X\to X$ and $\Bl_{Z_2} X\to X$, that is, the derived fibred product $\Bl_{Z_1} X\fibprod^\bR_X \Bl_{Z_2} X$ is a classical scheme.
The converse is not true: $\Bl_{\{0\}} \A^2\to \A^2$ is Tor-independent along itself but $\{0\}\to \A^2$ is not Tor-independent along itself.
This is a low-dimensional phenomenon though: $\Bl_{\{0\}} \A^n\to \A^n$ is not Tor-independent along itself for $n\geq 3$.
\end{exam}

\ssec{Local regular immersions}\label{ssec:loc-reg-emb}

\sssec{}
Let $i:Z\to X$ be a morphism of derived schemes.
We say that $i$ is \emph{unramified} if $i$ is locally of finite type and the relative cotangent complex $\sL_{Z/X}$ is $1$-connective (i.e., $\pi_{i}(\sL_{Z/X})=0$ for $i<1$).
This is equivalent to requiring that the underlying morphism of classical schemes is unramified (in the sense that $i_\cl : Z_\cl \to X_\cl$ is locally of finite type, and $\Omega_{Z_\cl/X_\cl} = 0$).
Equivalently, $i$ factors Zariski-locally on $Z$ as a closed immersion followed by an étale morphism~\cite[Cor.~18.4.7]{EGAIV4}.
In the category of derived algebraic spaces, such a factorization exists globally~\cite{rydh_embeddings-of-unramified}.

\begin{exam}\label{exam:finite+unramified}
A morphism of derived schemes is \emph{finite} if the underlying morphism of classical schemes is finite.
It follows from \cite[Cor.~18.4.7]{EGAIV4} that a morphism $i:Z\to X$ is finite and unramified if and only if, étale-locally on $X$, $i$ is a finite disjoint union of closed immersions.
\end{exam}

\sssec{}
A morphism $i : Z \to X$ is quasi-smooth and unramified if and only if it factors Zariski-locally on $Z$ as a quasi-smooth closed immersion followed by an étale morphism.
Alternatively, we have the following characterization, analogous to that of quasi-smooth closed immersions (\propref{prop:regular immersions and cotangent}).

\begin{prop}
Let $i:Z\to X$ be a morphism of derived schemes.
Then $i$ is quasi-smooth and unramified if and only if it is locally of finite presentation and the shifted cotangent complex $\sL_{Z/X}[-1]$ is a locally free $\sO_Z$-module of finite rank.
\end{prop}

\begin{proof}
We note that both conditions imply that $i$ is unramified.
Locally on $Z$, we can thus find a factorization $Z\hook X'\to X$ with $X'\to X$ \'etale and $Z\hook X'$ a closed immersion.
Then there is a canonical isomorphism $\sL_{Z/X} \simeq \sL_{Z/X'}$, so the result follows from \propref{prop:regular immersions and cotangent}.
\end{proof}

In analogy with the case of quasi-smooth immersions, we let $\sN_{Z/X} = \sL_{Z/X}[-1]$ and take this as the definition of the \emph{conormal sheaf} of a quasi-smooth unramified morphism.
The \emph{virtual codimension} of $Z\to X$ is the rank of $\sN_{Z/X}$.

\sssec{}
Let $X$ be a derived scheme.
A \emph{local virtual Cartier divisor} on $X$ is a derived scheme $D$ together with a quasi-smooth unramified morphism $i_D : D \to X$ of virtual codimension $1$.
Thus Zariski-locally on $D$ and étale-locally on $X$, a local virtual Cartier divisor is a virtual Cartier divisor.

\ssec{Derived blow-ups in unramified centres}

\sssec{}
Let $i : Z \to X$ be a quasi-smooth finite unramified morphism of derived schemes.
For any derived scheme $S$ and morphism $f : S \to X$, a \emph{virtual Cartier divisor on $S$ lying over} $(X,Z)$ is the datum of a commutative square
  \begin{equation*}
    \begin{tikzcd}
      D \ar{d}{g}\ar{r}{i_D}
        & S \ar{d}{f}
      \\
      Z \ar{r}{i}
        & X
    \end{tikzcd}
  \end{equation*}
satisfying the following conditions:

\begin{enumerate}
\item 
The morphism $i_D : D \to S$ exhibits $D$ as a local virtual Cartier divisor on $S$.

\item
The underlying square of classical schemes is cartesian.

\item
The canonical morphism
  \begin{equation*}
    g^*\sN_{Z/X} \to \sN_{D/S}
  \end{equation*}
is surjective (on $\pi_0$).
\end{enumerate}
Note that the second condition implies that $D\to S$ is finite.
It also implies that if $Z\to X$ is a closed immersion, then so is $D\to S$.
Thus in that case, the space of virtual Cartier divisors on $S$ over $(X,Z)$ and the space of local virtual Cartier divisors on $S$ over $(X,Z)$ coincide.

\sssec{}
As in \sssecref{sssec:Bl_Z/X}, we obtain a presheaf of spaces
  \begin{equation*}
    \Bl_{Z}X : (\DSch_{/X})^\op \to \Spc
  \end{equation*}
where $\Bl_{Z}X(S\to X)$ is the space of local virtual Cartier divisors over $(X,Z)$.
Note that if $Z\to X$ is a quasi-smooth closed immersion, then this definition agrees with the one in \sssecref{sssec:Bl_Z/X}.

\begin{thm}\leavevmode\label{thm:multiple}
Let $Z\to X$ be a quasi-smooth finite unramified morphism of derived schemes.
\begin{thmlist}
\item\label{thm:multiple/representable}
The derived stack $\Bl_{Z}X$ is (representable by) a derived algebraic space.

\item\label{thm:multiple/stable under base change}
The construction $\Bl_{Z} X \to X$ commutes with arbitrary derived base change, that is, $\left(\Bl_{Z} X X\right)\times^{\bR}_X X'=\Bl_{Z\times^{\bR}_X X'} X'$ for every morphism $X'\to X$ of derived schemes.

\item\label{thm:multiple/exceptional}
The diagonal $\Delta : Z\to Z\fibprod^\bR_X Z$ is an open and closed immersion.
Let $W=Z\fibprod^\bR_X Z\smallsetminus Z$ be its complement.
There is a canonical finite unramified morphism $\P_Z(\sN_{Z/X})\fibprod^\bR_Z \Bl_{W} Z \to \Bl_{Z} X$ which is the universal local virtual Cartier divisor lying over $(X,Z)$.

\item\label{thm:multiple/proper quasi-smooth}
The structural morphism $\pi_{Z/X} : \Bl_{Z} X \to X$ is proper and quasi-smooth.

\item\label{thm:multiple/disjoint}
When $Z=Z_1\amalg Z_2$, then $\Bl_{Z} X=\Bl_{Z_1} X\fibprod^\bR_X \Bl_{Z_2} X$.

\item\label{thm:multiple/etale}
If $g:X\to Y$ is an \'etale morphism such that $Z\to X\to Y$ is finite, then $\Bl_{Z} Y$ is the Weil restriction $g_*\Bl_{Z} X$.
\end{thmlist}
\end{thm}

\ssec{Proof of Theorem~\ref{thm:multiple}}
\sssec{Proof of~\itemref{thm:multiple/stable under base change}}
This is proven exactly as \thmref{thm:all}\itemref{thm:all/stable under base change}, see \sssecref{sssec:proof of stable under base change}.

\sssec{Proof of~\itemref{thm:multiple/disjoint}}

If $Z=Z_1\amalg Z_2$ is a disjoint union and $D$ is a local virtual Cartier divisor on $T$ lying over $(X,Z)$, then $D=D_1\amalg D_2$ where $D_i=g^{-1}(Z_i)$ is a local virtual Cartier divisor lying over $(X,Z_i)$. The result follows.

\sssec{Proof of~\itemref{thm:multiple/representable} and \itemref{thm:multiple/proper quasi-smooth}}

For a morphism $p:X'\to X$, let $i':Z'\to X'$ denote the derived base change of $i$ along $p$.
It is enough to show that $\Bl_{Z} X\times_X X'=\Bl_{Z'} X'$ is a derived scheme which is proper and quasi-smooth over $X'$ for some \'etale surjection $p$.
Since $i$ is finite, we can find an \'etale surjection $p$ such that $Z'\to X'$ is a disjoint union of closed immersions $i'_k:Z'_k\to X'_k$, see \sssecref{exam:finite+unramified}.
From~\itemref{thm:multiple/disjoint} and \thmref{thm:all}, we see that $\Bl_{Z'} X'=\Bl_{Z'_1} X'\fibprod^\bR_{X'} \dots \fibprod^\bR_{X'} \Bl_{Z'_n} X'$ is a derived scheme, quasi-smooth and proper over $X'$.

\sssec{Proof of~\itemref{thm:multiple/exceptional}}
We pick a smooth presentation as before.
It is enough to prove that the universal local virtual Cartier divisor of $\Bl_{Z'} X'$ is as stated.
The universal local virtual Cartier divisor of $\Bl_{Z'} X'$ over $Z'_k$ is the derived fibred product over $X'$ of $\P_{Z'_k}(\sN_{Z'_k/X'})$ and the $\Bl_{Z'_\ell} X'$ for all $\ell\neq k$.
This fibred product is isomorphic to the derived fibred product over $Z'_k$ of $\P_{Z'_k}(\sN_{Z'_k/X'})$ and the $\Bl_{Z'_\ell\fibprod^\bR_{X'} Z'_k/Z'_k}$ for all $\ell\neq k$.
The latter is $\Bl_{W'} Z'_k$ where $W'=\coprod_{\ell\neq k} \left(Z'_\ell\fibprod^\bR_{X'} Z'_k\right)=W\fibprod^\bR_X X'$.
The universal local virtual Cartier divisor of $\Bl_{Z} X$ is thus the derived fibred product over $Z$ of $\P_Z(\sN_{Z/X})$ and $\Bl_{W} Z$.

\sssec{Proof of~\itemref{thm:multiple/etale}}

Let $D$ be a local virtual Cartier divisor on $T$ over $(D,Y)$.
Then $D\to T$ factors as $D\to X\times_Y T\to T$.
Since $g : X\to Y$ is \'etale, $D\to X\times_Y T$ is also a local virtual Cartier divisor.
This gives a map $\Bl_{Z} Y\to g_*\Bl_{Z} X$.
Conversely, if $D$ is a local virtual Cartier divisor on $X\times_Y T$ over $(D,X)$, then the composite $D\to X\times_Y T\to T$ is a local virtual Cartier divisor over $(D,Y)$.
This gives a map $g_*\Bl_{Z} X\to \Bl_{Z} Y$ which is inverse to the previous map.


\bibliographystyle{alphamod}

{\small
\bibliography{references}
}

\end{document}